\documentstyle[12pt,amscd]{amsart}

\theoremstyle{plain} 

\newtheorem{ThNum}{Theorem}

\newtheorem{LemNum}{Lemma}

\newcommand{\hyp}{{\bold H}}

\newcommand{\real}{{\bold R}}
\newcommand{\complex}{{\bold C}}
\newcommand{\projective}{{\bold P}}

\newcommand{\PGL}{\operatorname{PGL}}

\newcommand{\Area}{\operatorname{Area}}

\begin{document}

\title[Configuration spaces and hyperbolic Dehn fillings, II]{%
	Configuration spaces of points on the circle and 
	hyperbolic Dehn fillings, II} 

\author[Y. Yamashita]{%
	Yasushi Yamashita} 
\address{%
	Department of Information and Computer Sciences \\
	Nara Women's University \\ 
	Kita-Uoya Nishimachi \\ 
	Nara 630-8506, Japan} 
\email{%
	yamasita@@ics.nara-wu.ac.jp} 
\author[H. Nishi]{%
	Haruko Nishi} 
\address{%
	Department of Mathematics \\
	Kyushu University \\ 
	33, Fukuoka 812-8581 Japan} 
\email{% 
	nishi@@math.kyushu-u.ac.jp} 
\author[S. Kojima]{%
	Sadayoshi Kojima} 
\address{%
	Department of Mathematical and Computing Sciences \\
	Tokyo Institute of Technology \\
	Ohokayama, Meguro \\    
	Tokyo 152-8552 Japan} 
\email{% 
	sadayosi@@is.titech.ac.jp}  

\date{% 
	Version 1.0 (November 26, 1998)} 

%\maketitle

%%%
%%%
%%%
\begin{abstract}
In our previous paper, we discussed the hyperbolization of the
configuration space of $n\,(\geq 5)$  marked points
with weights in the projective line up to projective transformations.
A variation of the weights induces a deformation.
It was shown that this correspondence of 
the set of the weights to the Teichm\"uller space when  $n = 5$  and
to the Dehn filling space when  $n= 6$
is locally one-to-one near the equal weight.
In this paper,
we establish its global injectivity.
\end{abstract}
\maketitle

%%%
%%%
%%%

\section{Introduction}
In \cite{KojimaNishiYamashita},
we have shown that the configuration space
of $n \, (\geq 5)$  marked points with weights 
on the real projective line up to projective transformations
admits a natural hyperbolization
so that the result becomes a hyperbolic cone-manifold
of dimension  $n-3$.
In brief,
we identify each component of the configuration space
with the space of similarity classes of convex $n$-gons with fixed
external angles  
in the complex plane via Schwarz-Christoffel map.
Then there is a beautiful way by Thurston to hyperbolize
such a space as an interior of some hyperbolic polyhedron
(see \cite{BavardGhys, KojimaYamashita}).
Each point on the boundary can be encoded by
an appropriately degenerate configuration.
Pasting them together along the same degenerate
configurations,
we obtained the resultant cone-manifold. 
Kapovich and Millson discussed the same hyperbolization 
via their duality in \cite{KapovichMillson}.  

The external angles can be given in fact at
your choice and we regarded them as the weight.
A variation of the weights induces a deformation.
We restricted the set of possible external
angles so that the $n$-gons are convex
and can be represented as the inner polygon of
the star shaped $n$-gons for any marking.
More concretely, we set 
\begin{equation*}
        \Theta_n=\{(\theta_1, \ldots, \theta_n)\,| \,
        \sum_{i=1}^n \theta_i= 2\pi, \,
        \theta_i>0, \, \theta_i + \theta_j < \pi \text { for any }i \ne j\}. 
\end{equation*}
This is equivalent to say that
the number of faces appeared in Thurston's polyhedralization
is constant,  that is  $n$.
Under this assumption,
the topology of a deformation will be almost constant.

In  \cite{KojimaNishiYamashita}, 
we discussed local behavior of the deformations
appeared in our setting near
the equal weight.   
When  $n = 5$, the deformations are topologically
a connected sum of five copies of the real projective 
space,  ${\#}^5 \real \projective^2$.
The assignment of the hyperbolic structure of a deformation
to each weight was shown to be a local embedding at the equal weight.
When  $n = 6$  and with the equal weight,
the result of hyperbolization is a $3$-dimensional hyperbolic
manifold of finite volume with ten cusps, which
we denoted by  $\overline{X_6}$.
Any deformation induced by a variation of the weights
can be regarded as some Dehn filled resultant of  $\overline{X_6}$.
The assignment of the deformation to each weight 
was also shown to be a local embedding at the equal weight.

In this paper,
we prove the global injectivity of the above assignment.
Namely, we show that  $\Theta_n$ is mapped by the above assignment
injectively to the deformation space in Theorem 1 when  $n = 5$,
and in Theorem 2  when  $n = 6$.
The local injectivity in \cite{KojimaNishiYamashita}
is proven by computing the derivative of the map at the equal weight.
The proof of the global injectivity we present here
is based on rather 
geometric observation for variation of polygons developed in
\cite{KojimaYamashita, AharaYamada}, and 
independent of the argument in \cite{KojimaNishiYamashita}.

We review some of materials in \cite{KojimaNishiYamashita} to set up
the notations in the next section,
and prove the theorems in the sections after.

\section{Preliminaries}

We here very briefly recall the hyperbolization 
in \cite{KojimaNishiYamashita}.  

The configuration space of $n$ 
marked points in the real projective line $\real \projective^1$ is, by
definition, the quotient of  $(\real \projective^1)^n$  minus 
the big diagonal set 
by the diagonal action of the projective
linear group $\PGL(2,\real)$.
It has $(n-1)!/2$ connected components,  
each of which is homeomorphic to a cell of dimension  $n-3$.  
Reading off the markings of 
the points in counterclockwise order, 
each component can be labeled by a circular permutation $p=\langle i_1 i_2
\ldots i_n \rangle$ of $n$ numbers $1$ to $n$ up to reversing the order.  

Fix an element $\theta = (\theta_1, \ldots, \theta_n)$ of $\Theta_n$.
Then there is a one-to-one Schwarz-Christoffel 
correspondence between 
the configuration space and the set of
similarity classes $X_{n,\theta}$ of marked $n$-gons in the complex plane
with external angles $\{\theta_1, \cdots, \theta_n\}$ compatible with markings
(see Lemma 1 in \cite{KojimaNishiYamashita}).
Fix a label $p=\langle i_1 i_2 \ldots i_n \rangle$.
Then the component of $X_{n,\theta}$ with label $p$ can be identified with
the subset of the space of all the congruence classes of Euclidean $n$-gons 
with external angles $\theta_{i_1}, \cdots, \theta_{i_n}$ cyclically
which consists of the ones with area 1.

Let $x_i$ denote the edge of 
an $n$-gon starting from the vertex with angle  $\theta_i$  in 
counterclockwise order, and simultaneously its length.  
Then, we have 
$$
\sum_{j=1}^n x_{i_j} \exp (\sum_{k=1}^j \sqrt{-1} \theta_{i_k}) =0.
$$

Let ${\cal E}_{p,\theta}$ be the $(n-2)$-dimensional vector space
satisfying the above constrain.
Then the space of congruence classes of $n$-gons is identified with 
the polyhedral 
cone ${\cal E}_{p,\theta}\cap\,\bigcap_{j=1}^n \{x_{i_j} >0\}$
in ${\cal E}_{p,\theta}$.  
The area determines a quadratic 
form $\Area$ of signature  $(1, n-3)$ on ${\cal E}_{p,\theta}$
(see Lemma 2 in \cite{KojimaNishiYamashita}).
Thus ${\cal E}_{p,\theta}$ together with $\Area$ becomes a Minkowski
space and $\Area^{-1}(1)$ is the hyperbolic
space in dimension $n-3$.
Therefore the space of similarity classes of $n$-gons
with fixed external 
angles $\theta_{i_1}, \cdots, \theta_{i_n}$ lies 
in the hyperbolic space bounded by the hyperbolic hyperplanes
$\Area^{-1}\cap \{x_{i_j}=0\}$ for $j = 1,\cdots, n$.
We denote by $\Delta_{p,\theta}$ such a hyperbolic polyhedron.
Then the conditions for $\theta$ as in the definition of $\Theta_n$
ensures us that the hyperbolic polyhedron $\Delta_{p,\theta}$
has exactly $n$ facets.  

Let us denote by $(i_1 i_2)i_3 \ldots i_n$ or simply $(i_1 i_2)$ the face of
$\Delta_{p,\theta}$ represented 
by $\Delta_{p,\theta} \cap \{x_{i_1}=0\}$ since it
corresponds to the degenerate configurations where 
the points marked by  $i_1$ and $i_2$
collide. Similarly we use $(i_1i_2)(i_3 i_4)i_5 \ldots i_n$ or 
$(i_1 i_2 i_3) i_4 \ldots i_n$, etc.
to represent the codimension two faces of $\Delta_{p,\theta}$.

Now gluing $(n-1)!/2$ hyperbolic polyhedra
$\Delta_{p, \theta}$ for all labels $p$ along the faces which 
represent the same degenerate configurations, 
we obtain $\overline{X_{n,\theta}}$ in which 
$X_{n, \theta}$ lies as an open dense subset.

\section{ The case $n=5$}

When $n=5$, $\Delta_{p, \theta}$ is a hyperbolic right pentagon
where the edges are labeled by $(i_1i_2)$, $(i_3i_4)$, $(i_5i_1)$, $(i_2i_3)$,
$(i_4i_5)$ cyclically. 
For the equal weight 
$\theta_0 =(2\pi/5, \cdots, 2\pi/5)$, $\Delta_{p, \theta_0}$ is a hyperbolic
regular pentagon, i.e, with all edges having equal lengths.

For any $\theta \in \Theta_5$, 
the space $\overline{X_{5,\theta}}$ is  a closed hyperbolic surface 
homeomorphic to ${\#}^5 \real \projective^2$.
The space of
hyperbolic structures, 
which we called a Teichm\"uller space and denoted 
by  ${\cal T}({\#}^5 \real \projective^2)$, 
is parameterized by the lengths and twisting amounts for the 2-sided ones
of a maximal family of mutually disjoint nonparallel 
simple closed curves on ${\#}^5 \real \projective^2$.
It is homeomorphic to $\real^9$.
Thus we have a map $\Phi_5: \Theta_5 \to {\cal T}({\#}^5 \real
\projective^2)$ assigning to each $\theta$ the marked 
hyperbolic structure of  $\overline{X_{5,\theta}}$.

On the other hand, our surface $\overline{X_{5,\theta}}$ has a 
geometric cell decomposition into twelve hyperbolic right pentagons. 
The edges of the pentagons form geodesics and they are uniquely placed on
$\overline{X_{5,\theta}}$ since the geodesic representative of 
simple closed curves
within their homotopy class is unique. 
Hence such a cell decomposition is uniquely
determined by the hyperbolic structure. In particular, the shapes of such pentagons
are invariants of the hyperbolic structure on
$\overline{X_{5,\theta}}$. Thus if we denote by ${\cal T}$ the space of all
hyperbolic right pentagons, we have a map from the subset $\Phi_5(\Theta_5)$
of ${\cal T}({\#}^5 \real \projective^2)$ to the direct 
product of twelve  ${\cal T}$'s  by listing 
hyperbolic structures of  $\Delta_{p, \theta}$'s. 
Here is a part of our surface (Fig.~\ref{Fig:1}).
\begin{figure}[ht]
\begin{center}
\setlength{\unitlength}{0.00083300in}%
\begin{picture}(4662,2724)(451,-4123)
\thicklines
\multiput(1351,-3211)(50,0){4}{.}
\multiput(1501,-3211)(40,-20){16}{.}
\multiput(2101,-3511)(0,-50){4}{.}
\multiput(3901,-1861)(0,-50){4}{.}
\multiput(3901,-2011)(40,-20){16}{.}
\multiput(4501,-2311)(50,0){4}{.}
\multiput(1351,-2761)(45,0){74}{.}
\multiput(2101,-1861)(0,-50){4}{.}
\multiput(2101,-2011)(-40,-20){16}{.}
\multiput(1501,-2311)(-50,0){4}{.}
\multiput(4651,-3211)(-50,0){4}{.}
\multiput(4501,-3211)(-40,-20){16}{.}
\multiput(3901,-3511)(0,-50){4}{.}
\multiput(3001,-3661)(0,45){41}{.}
\multiput(1501,-3211)(-30,-30){6}{.}
\multiput(2101,-3511)(-30,-30){6}{.}
\multiput(4051,-1861)(-30,-30){6}{.}
\multiput(4651,-2161)(-30,-30){6}{.}
\put(901,-1711){\line( 1,-1){600}}
\put(1501,-2311){\line( 0,-1){900}}
\put(2101,-3511){\line( 1, 0){1800}}
\put(3901,-3511){\line( 1,-1){600}}
\put(1501,-1411){\line( 1,-1){600}}
\put(2101,-2011){\line( 1, 0){1800}}
\put(4501,-2311){\line( 0,-1){900}}
\put(4501,-3211){\line( 1,-1){600}}
\put(3076,-3736){213(45)}
\put(2251,-1936){124(35)}
\put(3076,-1936){214(35)}
\put(2176,-2461){12435}
\put(3451,-2461){21435}
\put(4576,-2611){143(52)}
\put(826,-2611){243(51)}
\put(826,-3061){234(51)}
\put(2176,-3211){12345}
\put(3451,-3211){21345}
\put(4201,-3661){23145}
\put(4951,-3586){314(52)}
\put(4576,-3061){134(52)}
\put(2251,-3736){123(45)}
\put(451,-2086){423(51)}
\put(1276,-1936){14235}
\put(1726,-1561){142(35)}
\put(3601,-4036){231(45)}
\end{picture}
\end{center}
\caption{}\label{Fig:1}
\end{figure}

A hyperbolic right pentagon is determined by 
lengths of two adjacent sides. 
The parameterization of the space ${\cal T}$ can be given as follows.  
Put a hyperbolic right pentagon
in the projective model on the unit disk 
so that it lies in the first quadrant and the preferred two adjacent
sides are on the axes.
Let $P$ and $Q$ be the Euclidean lengths of the two sides on the axes.
Then these are subject to the relation $P^2 + Q^2 >1$. Conversely
if we are given such $P<1$ and $Q<1$, we can construct a unique hyperbolic pentagon.
Thus ${\cal T}$ can be identified with the set
\begin{equation*} 
	\{(P, Q)\;|\; 0 < P, Q < 1, P^2+Q^2>1\}. 
\end{equation*} 

We now review from \cite{KojimaYamashita} how the lengths $P$ and $Q$ 
for the hyperbolic right pentagon $\Delta_{p,\theta}$ can be computed 
in a geometric way where $p=\langle i_1i_2i_3i_4i_5\rangle$.
First recall that the space ${\cal E}_{p, \theta}$ of congruence classes
of Euclidean (extended) pentagons with external angles $\theta_{i_1}, \ldots,
\theta_{i_5}$ cyclically for $\theta = (\theta_1, \ldots, \theta_5)$ is a
$3$-dimensional Lorentz space. Its coordinate $(x, u, v)$ can be given 
by setting 
\begin{equation*}
x = \sqrt{\Area T_0}, \quad
u = \sqrt{\Area T_1}, \quad 
v = \sqrt{\Area T_2},
\end{equation*}
where $T_0$ is a triangle obtained by completing the pentagon $H$
by extending the edges $x_{i_2}$, $x_{i_4}$ and $x_{i_5}$,
and $T_1$, $T_2$ the yielding triangles adjacent to the
edge $x_{i_1}$, $x_{i_3}$ respectively (see Fig.~\ref{Fig:2}).

\begin{figure}[ht]
\begin{center}
\setlength{\unitlength}{0.00083300in}%
\begin{picture}(3624,1224)(1189,-2173)
\thicklines
\put(1201,-2161){\line( 1, 0){3600}}
\put(4801,-2161){\line(-2, 1){2400}}
\put(2401,-961){\line(-1,-1){1200}}
\put(3001,-2161){\line( 2, 5){300}}
\put(2101,-1261){\line( 1,-3){300}}
\put(2926,-1081){$T_0$}
\put(3556,-1928){$T_2$}
\put(1778,-1905){$T_1$}
\put(2108,-2064){$\theta_{i_2}$}
\put(3150,-2063){$\theta_{i_3}$}
\put(2618,-1674){$H$}
\put(1981,-1607){$\theta_{i_1}$}
\put(3300,-1681){$\theta_{i_4}$}
\end{picture}
\end{center}
\caption{}\label{Fig:2}
\end{figure}

Then $\Delta_{p,\theta}$ is a hyperbolic right pentagon 
lying in the hyperbola $\Area ^{-1}(1)$ in ${\cal E}_{p,\theta}$.
Projecting it along the rays towards the origin into the unit disk of the plane
$\{x=1\}$, we represent it in the projective model with $(u,v)$-coordinate.
Then we see that $\Delta_{p,\theta}$ lies 
in the first quadrant of the $uv$-plane
with the edge $(i_1i_2)i_3i_4i_5$ lying on the $u$-axis and
$i_1i_2(i_3i_4)i_5$ on the $v$-axis.
Let $(P,0)$ be the Euclidean coordinates of the intersection of
the adjacent edges $(i_1i_2)i_3i_4i_5$ and $i_1i_2i_3(i_4i_5)$, and $(0,Q)$ 
the intersection of adjacent edges $i_1i_2(i_3i_4)i_5$ and
$(i_5i_1)i_2i_3i_4$.  
To compute $P$, consider the degenerate pentagon $H_P$ with
$x_{i_1}=x_{i_4}=0$ in $T_0$ so that $T_0$ consists of $H_P$ and $T_2^e$, and
for $Q$ consider the degenerate pentagon $H_Q$ with degenerate sides
$x_{i_5}=x_{i_3}=0$ in $T_0$ so that $T_0$ consists of $H_Q$ and $T_1^e$ (see Fig.
\ref{Fig:3}).
\begin{figure}[ht]
\begin{center}
\setlength{\unitlength}{0.00083300in}%
\begin{picture}(5424,1224)(289,-4573)
\thicklines
\put(301,-4561){\line( 1, 0){2400}}
\put(2701,-4561){\line(-1, 1){1200}}
\put(1501,-3361){\line(-1,-1){1200}}
\put(1025,-4551){\line( 2, 5){476}}
\put(3301,-4561){\line( 1, 0){2400}}
\put(5701,-4561){\line(-1, 1){1200}}
\put(4501,-3361){\line(-1,-1){1200}}
\put(4501,-3361){\line( 1,-3){397}}
\put(1576,-4261){$T_2^e$}
\put(2101,-3661){$T_0$}
\put(4126,-4261){$T_1^e$}
\put(5176,-3661){$T_0$}
\put(4951,-4261){$H_Q$}
\put(863,-4246){$H_P$}
\end{picture}
\end{center}
\caption{}\label{Fig:3}
\end{figure}
Then by the definition of the coordinates of ${\cal E}_{p,\theta}$ and 
recalling that we are in the projective model in the plane $x=1$, we obtain 
\begin{align*}
P^2 & = \frac{\Area T_2^e}{\Area T_0} 
      = \frac{\text{ base width of }T_2^e}
             {\text{ base width of }T_0}
\end{align*}
and
\begin{align*}
Q^2 & = \frac{\Area T_1^e}{\Area T_0}
      = \frac{\text{ base width of }T_1^e}
             {\text{ base width of }T_0}. 
\end{align*}

Therefore by normalizing $T_0$ so that its base width has length 1,
the quantities $P^2$ and $Q^2$ are depicted just on the base edge
of $T_0$ as in Fig.~\ref{Fig:4}.
\begin{figure}[ht]
\begin{center}
\setlength{\unitlength}{0.00083300in}%
\begin{picture}(3924,1811)(889,-2473)
\thicklines
\put(901,-2161){\vector( 1, 0){3900}}
\put(1201,-2461){\vector( 0, 1){1787}}
\put(1201,-2161){\line( 1, 1){1200}}
\put(2401,-961){\line( 3,-2){1800}}
\put(2401,-961){\line(-1,-4){300}}
\put(2401,-961){\line( 1,-2){600}}
\put(4126,-2386){$1$}
\put(1951,-2386){$1-P^2$}
\put(2851,-2386){$Q^2$}
\put(3226,-961){$T_0$}
\put(976,-2386){$0$}
\end{picture}
\end{center}
\caption{}\label{Fig:4}
\end{figure}

Now let $\Psi_p : \Theta_5 \to {\cal T}$ be the map 
sending $\theta \in \Theta_5$  to 
the hyperbolic right pentagon $\Delta_{p,\theta}$.
Then the following is shown in \cite{KojimaYamashita}.

\begin{LemNum}\label{Lem:5fibration}
The inverse image 
of a point in ${\cal T}$  by $\Psi_p$  is in one-to-one
correspondence with the upper halfplane $\hyp^2$.
\end{LemNum}

We review here the proof in \cite{KojimaYamashita} for 
reader's convenience.

\begin{pf}
Suppose $(P, Q) \in {\cal T}$.
For a point $w$ in $\hyp^2$, plot four points $0 < 1-P^2 < Q^2 < 1$ 
on the real axis of $\complex$ and join these points with $w$
(Fig.~\ref{Fig:5}). We denote these four edges by $e_0, e_P, e_Q, e_1$ as in the
figure. Slide $e_P$ and $e_Q$ to the points $1$ and $0$
respectively so that the lines are parallel to the originals and do
not intersect with each other.  Then these new two edges together with
$e_0$, $e_1$ and the real--axis form a pentagon with its external angles
$(\theta_{i_1}, \ldots, \theta_{i_5})$ as in Fig.~\ref{Fig:5}.
Then by the argument above, we see that $\theta=
(\theta_1, \ldots, \theta_5)$ satisfies $\Psi_p(\theta) =(P, Q)$.
Since the choice of $w$ is arbitrary in $\hyp ^2$, by checking that the
different $w$ gives different $\theta$ we get the conclusion.
\end{pf}
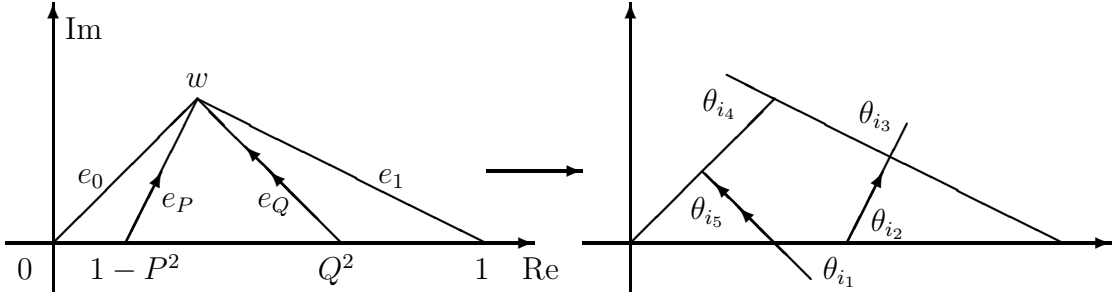
\begin{figure}[ht]
\begin{center}
\setlength{\unitlength}{0.00083300in}%
\begin{picture}(6924,1824)(289,-1273)
\put(2626,-586){$e_1$}
\thicklines
\put(601,-1261){\vector( 0, 1){1800}}
\put(601,-961){\line( 1, 1){900}}
\put(1501,-61){\line( 2,-1){1800}}
\put(3901,-961){\vector( 1, 0){3300}}
\put(4201,-1261){\vector( 0, 1){1800}}
\put(1501,-61){\line(-1,-2){450}}
\put(3301,-511){\vector( 1, 0){600}}
\put(4201,-961){\line( 1, 1){900}}
\put(6901,-961){\line(-2, 1){2100}}
\put(1051,-961){\vector( 1, 2){225}}
\put(2401,-961){\line(-1, 1){900}}
\put(5326,-1186){\line(-1, 1){675}}
\put(5551,-961){\vector( 1, 2){225}}
\put(5926,-211){\line(-1,-2){375}}
\put(2101,-661){\vector(-1, 1){300}}
\put(2251,-811){\vector(-1, 1){300}}
\put(5026,-886){\vector(-1, 1){300}}
\put(5176,-1036){\vector(-1, 1){300}}
\put(826,-1186){$1-P^2$}
\put(1426, 14){$w$}
\put(376,-1186){$0$}
\put(3226,-1186){$1$}
\put(676,314){Im}
\put(2251,-1186){$Q^2$}
\put(4576,-811){$\theta_{i_5}$}
\put(5401,-1186){$\theta_{i_1}$}
\put(5701,-886){$\theta_{i_2}$}
\put(4651,-136){$\theta_{i_4}$}
\put(3526,-1186){Re}
\put(5626,-211){$\theta_{i_3}$}
\put(751,-586){$e_0$}
\put(1276,-736){$e_P$}
\put(1876,-736){$e_Q$}
\put(301,-961){\vector( 1, 0){3300}}
\end{picture}
\caption{Map from $w$ to $\theta$}\label{Fig:5}
\end{center}
\end{figure}

Now consider the product map $\Psi_{\langle 12345\rangle} \times \Psi_{\langle 21435\rangle}: 
\Theta_5 \to {\cal T}\times {\cal T}$ which sends $\theta$ to 
$(\Delta_{\langle 12345\rangle,\theta},\Delta_{\langle 21435\rangle,\theta})$.
Then we claim that 

\begin{LemNum}\label{Lem:5injective}
The map $\Psi_{\langle 12345\rangle} \times \Psi_{\langle 21435\rangle}$ is injective.
\end{LemNum}

\begin{pf}
Let $((P_1, Q_1), (P_2, Q_2))$ be any element of the image of 
$\Psi_{\langle 12345\rangle} \times \Psi_{\langle 21435\rangle}$, that is, $P_1$, $Q_1$, $P_2$, $Q_2$ are the
lengths of the sides 
$(12)345$, $12(34)5$, $(21)435$, $21(43)5$, respectively.   Suppose that
$\theta = (\theta_1,\ldots,\theta_5)$ be any element in
the inverse image of $((P_1, Q_1), (P_2, Q_2))$ by 
$\Psi_{\langle 12345\rangle} \times \Psi_{\langle 21435\rangle}$. 
Then $\theta$ lies in the intersection 
$\Psi_{\langle 12345\rangle}^{-1}((P_1,Q_1)) \cap \Psi_{\langle 21435\rangle}^{-1}((P_2,Q_2))$. 
Let $w_1$ be the point in $\hyp^2 \cong \Psi_{\langle 12345\rangle}^{-1}((P_1,Q_1))$
corresponding to $\theta$ under the identification in Lemma~\ref{Lem:5fibration},
and similarly $w_2$ the point in $\hyp^2 \cong \Psi_{\langle 21435\rangle}^{-1}((P_2,Q_2))$ 
corresponding to $\theta$.
These can be obtained by the inverse procedure of Fig.~\ref{Fig:5}.
Then we see that 
the triangles $\triangle w_101$ and $\triangle w_201$ are congruent since both have
the external angles $\theta_5$, $\theta_1+\theta_2$, $\theta_3+\theta_4$.
Thus $w_1$ and $w_2$ must be identical in the upper half plane $\hyp^2$.
Let $w=w_1=w_2$. 

The angle at $P_1^2$ in the triangle $\triangle w0P_1^2$ is
$\theta_2$,  and 
the angle at $P_2^2$ in the triangle $\triangle w0P_2^2$ 
is $\theta_1$, so
$\triangle w 0P_1^2$ and $\triangle P_2^20w$ are similar.
(See Fig.~\ref{Fig:6}).
\begin{figure}[ht]
\begin{center}
\setlength{\unitlength}{0.00083300in}%
\begingroup\makeatletter\ifx\SetFigFont\undefined%
\gdef\SetFigFont#1#2#3#4#5{%
  \reset@font\fontsize{#1}{#2pt}%
  \fontfamily{#3}\fontseries{#4}\fontshape{#5}%
  \selectfont}%
\fi\endgroup%
\begin{picture}(6912,1767)(1,-1597)
%\put(6601,-886){$\theta_4+\theta_3$}
\put(6301,-886){$\theta_4+\theta_3$}
\thicklines
\put(3901,-961){\vector( 1, 0){3000}}
\put(601,-961){\line( 1, 1){900}}
\put(1501,-61){\line( 4,-3){1200}}
\put(1501,-61){\line( 1,-3){300}}
%\put(5101,-61){\line( 5,-3){1500}}
\put(5101,-61){\line( 4,-3){1200}}
\put(5101,-61){\line( 2,-3){600}}
\put(4201,-961){\line( 1, 1){900}}
\put(5101,-961){\line( 0, 1){900}}
\put(1201,-961){\line( 1, 3){300}}
\put(5626,-1186){$P_2^2$}
\put(1726,-1186){$P_1^2$}
\put(5476,-1561){(b)}
\put(1876,-1561){(a)}
\put(526,-1186){$0$}
\put(4126,-1186){$0$}
\put(6226,-1186){$1$}
\put(2626,-1186){$1$}
\put(1576,-886){$\theta_2$}
\put(5401,-886){$\theta_1$}
\put(5176,-886){$\theta_4$}
\put(1426, 14){$w_1$}
\put(4951, 14){$w_2$}
\put(1276,-886){$\theta_3$}
\put(  1,-886){$\theta_1+\theta_2$}
\put(2701,-886){$\theta_3+\theta_4$}
\put(3601,-886){$\theta_2+\theta_1$}
\put(976,-1186){$1-Q_1^2$}
\put(4801,-1186){$1-Q_2^2$}
\put(301,-961){\vector( 1, 0){2700}}
\end{picture}
\caption{}\label{Fig:6}
\end{center}
\end{figure}
Then we have
\[ |w| :  P_1^2 =  P_2^2 : |w|. \]
Therefore
\[ |w|^2 =  (P_1 P_2)^2. \]
This implies that $w$ must lie on the circle with center
$0$ and radius $P_1P_2$. In the same way, since the triangles
$\triangle w_1(1-Q_1^2)1$ and $\triangle (1-Q_2^2)1w_2$ are similar, we see that
$w$ must lie on the circle with center $1$ and radius $Q_1Q_2$. Since the
above two circles intersect at one point on the upper halfplane, the point $w$,
therefore $\theta$, is uniquely determined.
\end{pf}

By this lemma, we see that the composition
assigning $\theta \mapsto \overline{X_{5,\theta}}
\mapsto (\Delta_{\langle 12345\rangle, \theta}, 
\Delta_{\langle 21435\rangle, \theta})$ is injective, and 
hence we have proven the following.

\begin{ThNum}
The map $\Phi_5$ from $\Theta_5$ to ${\cal T}({\#}^5 \real \projective^2)$  
defined by the assignment to
$\theta$ the hyperbolic structure of $\overline{X_{5,\theta}}$ is
injective.
\end{ThNum}

\section{ The case $n=6$ }

When $n=6$, $\Delta_{p, \theta}$ is a hyperbolic hexahedron 
bounded by six faces $(i_1 i_2), (i_2i_3)\ldots, (i_6 i_1)$ with
having two vertices $(i_1 i_2)(i_3 i_4)(i_5 i_6)$
and $(i_2 i_3)(i_4 i_5)(i_6 i_1)$ at finite distance where
three faces $(i_k i_{k+1})$ are coming
and meeting orthogonally each other
for $k=1,3,5$ and $k=2,4,6$ respectively.
Moreover, the faces $(i_k i_{k+1})$ and $(i_{k+3} i_{k+4})$
always intersect orthogonally for $k=1,2,3$.
For the equal weight $\theta_0 =(2\pi/6, \cdots, 2\pi/6)$,
the dihedral angle between the faces $(i_k i_{k+1})$ and $(i_{k+1} i_{k+2})$
is zero for $k=1,\ldots, 6$ so that $\Delta_{p, \theta_0}$ has 
three ideal vertices.
For other weights $\theta \ne \theta_0$,
the faces $(i_k i_{k+1})$ and $(i_{k+1} i_{k+2})$ intersect
if and only if $\theta_{i_k}+ \theta_{i_{k+1}}+ \theta_{i_{k+2}} < \pi$.
Note that in such a case, $(i_{k+3} i_{k+4})$ and $(i_{k+4} i_{k+5})$
don't intersect since $\theta_{i_{k+3}}+ \theta_{i_{k+4}}+ \theta_{i_{k+5}} > \pi$.

The hyperbolized configuration space of the equal weight 
in this case is a complete hyperbolic 3-manifold $\overline{X_6}$
of finite volume with ten cusps.  
A perturbation gives rise to not a deformation 
in the usual sense, 
but a resultant of hyperbolic Dehn filling.  
The space of hyperbolic Dehn fillings of a complete 
hyperbolic 3-manifold can be locally identified with 
the space of representations 
of a fundamental group up to conjugacy.  
It has a structure of a smooth algebraic variety of 
complex dimension $= \text{the number of cusps}$, 
(\cite{CullerShalen},
\cite{ThurstonNote}, \cite{NeumannZagier}, \cite{HodgsonKerckhoff}).   
Hence in our case,  the space of Dehn
fillings ${\cal H}(\overline{X_6})$ is locally biholomorphic to  $\complex^{10}$.    

Here is a part of our 3--manifold $\overline{X_6}$
in the projective model (Fig.~\ref{Fig:7}). 
\begin{figure}[ht]
\begin{center}
\setlength{\unitlength}{0.00083300in}%
\begin{picture}(3698,2640)(2678,-3094)
\thicklines
\put(2901,-2462){\line( 2,-1){806}}
\put(3707,-2865){\line( 6, 1){2416.054}}
\put(4915,-1858){\line(-2, 5){402.483}}
\put(4915,-1858){\line(-6,-5){1208.164}}
\put(3304,-1657){\line( 1,-3){402.700}}
\put(3304,-1657){\line(-1,-2){402.600}}
\put(3304,-1657){\line( 3, 2){1207.846}}
\put(5721,-1456){\line( 2,-5){402.345}}
\put(5721,-1456){\line(-2, 1){1208.800}}
\put(4915,-1858){\line( 2,-1){1208}}
\put(2700,-2462){\vector( 1, 0){3625}}
\put(4512,-2865){\vector( 0, 1){2215}}
%% \put(3976,-1411){\vector( 1,-1){0}}
%% \put(3676,-1561){\vector( 2,-1){0}}
%% \put(4876,-1336){\vector(-1,-1){0}}
\put(5519,-1959){\vector(-2,-1){2014}}
\multiput(5721,-1456)(-33.56410,-50.34615){13}{.}
\multiput(4110,-1456)(52.52174,-26.26087){24}{.}
\multiput(4110,-1456)(-46.47163,-38.72636){27}{.}
\multiput(3601,-1036)(41.66667,-41.66667){9}{.}
\multiput(3076,-1261)(54.54545,-27.27273){11}{.}
\multiput(5326,-886)(-40.90909,-40.90909){11}{.}
\multiput(4512,-852)(-22.36590,-55.91475){37}{.}
\multiput(4512,-852)(42.38158,-42.38158){39}{.}
\multiput(4512,-852)(-42.38158,-42.38158){39}{.}
\multiput(2901,-2462)(58.94792,9.82465){42}{.}
\multiput(5318,-2060)(53.65333,-26.82667){16}{.}
\multiput(4110,-1456)(33.53846,50.30769){13}{.}
\multiput(4512,-852)(33.56410,-50.34616){25}{.}
\put(4778,-2986){123456}
\put(2678,-2836){124356}
\put(5603,-1261){213456}
\put(3278,-961){214356}
\put(3226,-3061){$u$}
\put(6376,-2536){$v$}
\put(4426,-586){$w$}
\put(5026,-2386){(23)}
%\put(4126,-2086){(61)}
%\put(3751,-1936){$A$}
\put(3951,-1836){(61)}
%\put(4966,-1786){(45)}
%\put(5326,-1561){$B$}
\put(5166,-1686){(45)}
%\put(5326,-811){$C (61)$}
\put(5326,-811){$(62)$}
%\put(2776,-1186){$D (45)$}
\put(2776,-1186){$(35)$}
\end{picture}
\caption{}\label{Fig:7}
\end{center}
\end{figure}
The figure illustrates the union of four hexahedra 
in the case of the equal weight.  
The vertices on the axes are located on the ideal boundary, but other 
four vertices are at finite distance.  
%One component of $X(6)$ corresponds to the
%intersection of the polyhedron with one orthant.  
Six planes $u=0$, $v=0$, $w=0$, $u+v=1$,
$v+w=1$ and $w+u=1$ correspond to the collision of 
$(12)$, $(34)$, $(56)$, $(23)$, $(45)$
and $(61)$ respectively. The origin corresponds to $(12)(34)(56)$.  

If we perturb the weight to $\theta$ so 
that $\theta_{i_k} + \theta_{i_{k+1}} +\theta_{i_{k+2}} < \pi$, 
then three vertices $i_{k}$, $i_{k+1}$,
$i_{k+2}$ of the hexagon can collide, and $(i_{k}i_{k+1})$ and
$(i_{k+1}i_{k+2})$ intersects to yield a new edge.
The dihedral angles about the old edges are $\pi/2$
(\cite{KojimaNishiYamashita}) so that the bounding planes form geodesic
surfaces.  
But the dihedral angle at the new edge depends on the weight $\theta$
and it creates a cone-type singularity in $\overline{X_{6,\theta}}$ so that 
$\overline{X_{6,\theta}}$ becomes a hyperbolic cone-manifold.
Notice that the complement of the cone-singularity is homeomorphic to
$\overline{X_6}$ but endowed with an 
incomplete hyperbolic structure.
Then $\overline{X_{6,\theta}}$ can be regarded as the
resultant of the hyperbolic Dehn filling of the complete 
hyperbolic 3-manifold $\overline{X_6}$, which 
lies in ${\cal H}(\overline{X_6})$
(for more detail, see \cite{KojimaNishiYamashita}). 
Thus we have a map $\Phi_6: \Theta_6 \to {\cal H}(\overline{X_6})$ assigning 
to each $\theta$ 
the hyperbolic cone-manifold $\overline{X_{6,\theta}}$.

There are 15 geodesic surfaces in $\overline{X_{6,\theta}}$ represented 
by a pair of
markings $(ij)$ where $i \ne j \in \{1, 2, \ldots, 6\}$.  
They are uniquely placed on
$\overline{X_{6,\theta}}$ since the geodesic representative of a
surface within their proper homotopy class is unique if any.  Hence a
geometric cell decomposition by such surfaces, which consists of 60
hexahedra, is uniquely determined and
the shapes of these hexahedra $\Delta_{p, \theta}$ are invariants of
the hyperbolic (cone) structure on $\overline{X_{6,\theta}}$.
Let ${\cal H}$ be the space of all hexahedra which have the
property as $\Delta_{p,\theta}$ described 
in the first paragraph of this section.
Then we have a map from 
the subset $\Phi_6(\Theta_6)$ of ${\cal H}$ to the direct 
product of sixty ${\cal H}$'s  by listing hyperbolic structures 
of  $\Delta_{p, \theta}$'s. 

The parameterization of the space ${\cal H}$ can be given as follows.
First let us make the definition of  ${\cal H}$  precise.
We use the notations of the faces or vertices of the elements in ${\cal H}$
as $\Delta_{p,\theta}$ where $p=\langle 123456\rangle$.
Then ${\cal H}$ is the set of all hyperbolic hexahedra having the following
properties; (i) it is bounded by six faces $(12),
(23)\ldots, (61)$,
(ii) the three faces $(k\,k+1)$ meet orthogonally each other
for $k=1,3,5$ and $k=2,4,6$ to make vertices $(12)(34)(56)$
and $(23)(45)(61)$ at finite distance, respectively,
(iii) the faces $(k\,k+1)$ and $(k+3\,k+4)$
always meet orthogonally for $k=1,2,3$.

We present a hyperbolic hexahedron in ${\cal H}$ in the projective model
located in the unit ball in the Euclidean 3-space with $uvw$-coordinates
so that the vertex $(12)(34)(56)$ is at the origin, 
the faces $(12)$, $(34)$, $(56)$ are lying on the $vw$, $wu$, $uv$-planes,
respectively. Since the other three faces $(45)$, $(61)$ and $(23)$ are 
orthogonal in a hyperbolic sense to the faces $(12)$, $(34)$ and $(56)$ respectively, 
they are orthogonal also in an Euclidean sense in this model there.
We denote the hyperplane on which the face $(ij)$ lies by the same symbol
$(ij)$, for simplicity.

Let $(P,0,0)$ be the Euclidean coordinates of the
intersection of the plane $(61)$ with $u$-axis,
$(0,Q,0)$ the intersection of $(23)$ with  $v$-axis,
$(0,0,R)$ the intersection of $(45)$ with  $w$-axis, respectively.
Then depending on the signs of $P-1$, $Q-1$, $R-1$, there are
several combinatorial types of hexahedra in ${\cal H}$ (see Fig.~\ref{Fig:8}.)
\begin{figure}[ht]
\begin{center}
\setlength{\unitlength}{0.00083300in}%
\begingroup\makeatletter\ifx\SetFigFont\undefined%
\gdef\SetFigFont#1#2#3#4#5{%
  \reset@font\fontsize{#1}{#2pt}%
  \fontfamily{#3}\fontseries{#4}\fontshape{#5}%
  \selectfont}%
\fi\endgroup%
\begin{picture}(7224,6457)(289,-6806)
\thicklines
\put(1501,-2161){\vector( 1, 0){1800}}
\put(1501,-2161){\vector(-1,-1){1200}}
\put(1501,-2161){\vector( 0, 1){1800}}
\put(901,-2761){\line( 1, 5){300}}
\put(1201,-1261){\line( 1, 1){300}}
\put(1501,-961){\line( 4,-3){1200}}
\put(2701,-1861){\line( 0,-1){300}}
\put(2701,-2161){\line(-5,-2){1500}}
\put(1201,-2761){\line(-1, 0){300}}
\put(1201,-2761){\line( 2, 3){600}}
\put(1801,-1861){\line( 1, 0){900}}
\put(1201,-1261){\line( 1,-1){600}}
\put(1501,-5461){\vector( 0, 1){1800}}
\put(1501,-5461){\vector(-1,-1){1200}}
\put(1501,-5461){\vector( 1, 0){2700}}
\put(1501,-4261){\line(-1,-1){300}}
\put(1201,-4561){\line(-1,-3){300}}
\put(901,-5461){\line( 0,-1){600}}
\put(901,-6061){\line( 5, 1){1500}}
\put(2401,-5761){\line( 1, 1){300}}
\put(2701,-5461){\line(-1, 1){1200}}
\put(1201,-4561){\line( 1,-2){300}}
\put(1501,-5161){\line(-2,-1){600}}
\put(1501,-5161){\line( 3,-2){900}}
\put(4801,-5461){\vector( 0, 1){1800}}
\put(4801,-5461){\vector(-1,-1){1200}}
\put(4801,-5461){\vector( 1, 0){2700}}
\put(4801,-4261){\line( 1, 0){300}}
\put(5101,-4261){\line( 3,-4){900}}
\put(6001,-5461){\line(-1,-1){300}}
\put(5701,-5761){\line(-5,-1){1500}}
\put(4201,-6061){\line( 0, 1){300}}
\put(4201,-5761){\line( 2, 5){600}}
\put(5101,-4261){\line(-1,-6){150}}
\put(4951,-5161){\line( 5,-4){750}}
\put(4951,-5161){\line(-5,-4){750}}
\put(4801,-2161){\vector( 0, 1){1800}}
\put(4801,-2161){\vector(-1,-1){1200}}
\put(4801,-2161){\vector( 1, 0){1800}}
\put(4801,-961){\line( 4,-3){1200}}
\put(6001,-1861){\line( 0,-1){300}}
\put(6001,-2161){\line(-3,-1){1800}}
\put(4201,-2761){\line( 0, 1){600}}
\put(4801,-961){\line(-1,-1){300}}
\put(4501,-1261){\line(-1,-3){300}}
\put(4501,-1261){\line( 1,-1){600}}
\put(5101,-1861){\line(-3,-1){900}}
\put(5101,-1861){\line( 1, 0){900}}
\put(751,-2761){\line( 0,-1){300}}
\put(4051,-6061){\line( 0,-1){300}}
\put(751,-6061){\line( 0,-1){300}}
\put(3001,-5311){\line( 0,-1){300}}
\put(6301,-5311){\line( 0,-1){300}}
\put(6301,-2011){\line( 0,-1){300}}
\put(3001,-2011){\line( 0,-1){300}}
\put(4051,-2761){\line( 0,-1){300}}
\put(1351,-661){\line( 1, 0){300}}
\put(4651,-661){\line( 1, 0){300}}
\put(4651,-3961){\line( 1, 0){300}}
\put(1351,-3961){\line( 1, 0){300}}
\put(376,-3511){$u$}
\put(3751,-3436){$u$}
\put(3226,-2386){$v$}
\put(6526,-2386){$v$}
\put(3751,-5686){$v$}
\put(7276,-5686){$v$}
\put(901,-6211){$p$}
\put(3601,-5386){$Q$}
\put(6976,-5386){$Q$}
\put(1051,-3511){(a) $P<1, Q<1, R<1$}
\put(4501,-3511){(b) $P>1, Q<1, R<1$}
\put(4501,-6736){(d) $P>1, Q>1, R>1$}
\put(1051,-6736){(c) $P>1, Q>1, R<1$}
\put(1651,-2386){(23)}
\put(1726,-1561){(45)}
\put(1126,-2086){(61)}
\put(601,-1411){(34)}
\put(2476,-1336){(12)}
\put(1801,-2986){(56)}
\put(5881,-2356){$q$}
\put(2603,-2356){$q$}
\put(4239,-2918){$p$}
\put(924,-2949){$p$}
\put(4224,-6226){$p$}
\put(586,-6556){$P$}
\put(399,-6773){$u$}
\put(3721,-6766){$u$}
\put(1599,-3818){$w$}
\put(4913,-3819){$w$}
\put(4891,-496){$w$}
\put(1599,-526){$w$}
\put(4614,-4284){$r$}
\put(6039,-5394){$q$}
\put(2708,-5386){$q$}
\put(3879,-6571){$P$}
\put(4569,-3804){$R$}
\put(3894,-3256){$P$}
\put(4898,-961){$r$}
\put(1576,-946){$r$}
\put(1583,-4262){$r$}
\put(1726,-719){$1$}
\put(2948,-2512){$1$}
\put(683,-3269){$1$}
\multiput(901,-5461)(-18.75000,-56.25000){17}{\makebox(6.6667,10.0000){\SetFigFont{10}{12}{\rmdefault}{\mddefault}{\updefault}.}}
\multiput(2401,-5761)(60.00000,12.00000){26}{\makebox(6.6667,10.0000){\SetFigFont{10}{12}{\rmdefault}{\mddefault}{\updefault}.}}
\multiput(4201,-2161)(-18.75000,-56.25000){17}{\makebox(6.6667,10.0000){\SetFigFont{10}{12}{\rmdefault}{\mddefault}{\updefault}.}}
\multiput(5101,-4261)(-33.33333,50.00000){10}{\makebox(6.6667,10.0000){\SetFigFont{10}{12}{\rmdefault}{\mddefault}{\updefault}.}}
\multiput(5701,-5761)(60.00000,12.00000){26}{\makebox(6.6667,10.0000){\SetFigFont{10}{12}{\rmdefault}{\mddefault}{\updefault}.}}
\multiput(4201,-5761)(-27.27273,-54.54545){12}{\makebox(6.6667,10.0000){\SetFigFont{10}{12}{\rmdefault}{\mddefault}{\updefault}.}}
\multiput(1801,-2766)(-42.85714,42.85714){7}{\makebox(6.6667,10.0000){\SetFigFont{10}{12}{\rmdefault}{\mddefault}{\updefault}.}}
\multiput(2393,-1424)(-42.71429,-42.71429){7}{\makebox(6.6667,10.0000){\SetFigFont{10}{12}{\rmdefault}{\mddefault}{\updefault}.}}
\multiput(976,-1409)(42.85714,-42.85714){7}{\makebox(6.6667,10.0000){\SetFigFont{10}{12}{\rmdefault}{\mddefault}{\updefault}.}}
\end{picture}
\end{center}
\caption{}\label{Fig:8}
\end{figure}
To describe each type, it is enough to look at the cases $(a)$ $P<1, Q<1,
R<1$ and $(b)$ $P>1, Q<1, R<1$.
Let $p =\text{ mim }\{P, 1/P\}$, $q=\text{ mim }\{Q, 1/Q\}$,
$r=\text{ mim }\{R, 1/R\}$. 

Case $(a)$. $P<1, Q<1, R<1$;

\noindent
In this case, since the face $(61)$ intersects with $(56)$,
the face $(23)$ with $(12)$, the face $(45)$ with $(34)$ within the unit ball in the
$uvw$-space respectively, all the faces of $\Delta_{p, \theta}$ are
quadrilaterals with three right angles. (See Fig.~\ref{Fig:8}$(a)$.)
Then $p =P$, $q =Q$, $r =R$, 
and the points $(p,0,0)$, $(0,q,0)$, $(0,0,r)$ are the vertices
of the hexahedron.

Case $(b)$. $P>1, Q<1, R<1$;

\noindent
In this case, the plane $(61)$ does not intersect with
$(56)$ inside the unit ball, while the faces $(23)$ and 
$(34)$ intersect. (See Fig.~\ref{Fig:8}$(b)$).
Then the face $(34)$ is a hyperbolic right pentagon,
whose side on the $u$-axis has the Euclidean length $p=1/P$.
This is because when restricted to the hyperbolic 2-plane $\{v=0\}$,
the point $(P,0,0)$ is the intersection of the two
lines extending the edges $(61)\cap(34)$ and $(34)\cap(56)$ respectively,
both of which are orthogonal to the edge $(23)\cap(34)$.
In other words, the point $(P,0,0)$ is the dual of the
line containing the edge $(23)\cap(34)$ (see Fig~\ref{Fig:9}).
\setlength{\unitlength}{0.00083300in}%
\begin{figure}[ht]
\begin{center}
\setlength{\unitlength}{0.00083300in}%
\begin{picture}(3038,2407)(1853,-1856)
\thicklines
\put(4801,-1561){\vector(-1, 0){2700}}
\put(4801,-1561){\vector( 0, 1){2100}}
\put(4801,239){\line(-5,-4){2250}}
\put(4726,-61){\line( 1, 0){150}}
\put(3301,-1486){\line( 0,-1){150}}
\put(4651,-1561){\line( 0, 1){150}}
\put(4651,-1411){\line( 1, 0){150}}
\put(4801,-436){\line(-1, 0){150}}
\put(4651,-436){\line( 0, 1){150}}
\put(3676,-1411){\line( 1, 0){150}}
\put(3826,-1411){\line( 0,-1){150}}
\put(3676,-1561){\line( 0, 1){900}}
\put(4801,-286){\line(-1, 0){650}}
\put(3676,-736){\line( 1, 0){150}}
\put(3826,-736){\line( 0, 1){190}}
\put(4201,-286){\line( 0,-1){150}}
\put(4201,-436){\line(-1, 0){225}}
\put(4891,-338){$r$}
\put(4516,337){$w$}
\put(1853,-1613){$u$}
\put(4891,-113){$1$}
\put(3226,-1823){$1$}
\put(2498,-1809){$P$}
\put(3481,-1823){$1/P$}
\end{picture}
\end{center}
\caption{}\label{Fig:9}
\end{figure}
The shape of the faces $(23)$ is also a hyperbolic pentagon,
the faces $(12)$ and $(45)$ are hyperbolic quadrilaterals with three
right angles, and the faces $(61)$ and $(56)$ are hyperbolic right triangles.

Note that when $P=1, Q<1, R<1$, the hexahedra has an ideal vertex at $(1,0,0)$.
This means that the faces $(61)$ and $(56)$ do not meet, but
they are tangent at $(1,0,0)$, and also $(23)$ and $(34)$ are
tangent at $(1,0,0)$.
So the face $(61)$ and the face $(56)$ are hyperbolic right triangles with one ideal vertex, 
the face $(23)$ and the face $(34)$ quadrilaterals with one ideal vertex and three right angles.
The faces $(12)$ and $(45)$ are as in the case $(a)$ or $(b)$.

Now conversely when three numbers $P$, $Q$, $R >0$ are given, 
by taking the inverse procedure above,
we can construct a unique hyperbolic hexahedron in ${\cal H}$.
That is, its combinatorial type is determined by the signs of $P-1$,
$Q-1$, $R-1$, and the lengths of the edges
on the $u$, $v$, $w$-axes are $p$, $q$, $r$,
respectively.
Since each face is
either a hyperbolic right pentagon, a quadrilateral with three right 
angles (possibly with one ideal vertices),
or a hyperbolic right triangle (possibly with ideal vertices), 
all the other lengths of the hexahedron are uniquely determined
and the vertex where $(61)$, $(23)$ and $(45)$ meet is 
at finite distance.
Therefore ${\cal H}$ can be identified with the set 
\begin{equation*} 
	\{P, Q, R \, | \, P, Q, R  >0\}. 
\end{equation*} 

Following \cite{AharaYamada}, we shall show 
a geometric way of computing the parameters $P$, $Q$, $R$
for the hyperbolic hexahedron $\Delta_{p,\theta}$
as in the case $n=5$.  

Recall that for $p=\langle i_1i_2i_3i_4i_5i_6\rangle$ the 
space ${\cal E}_{p,\theta}$ of 
hexagons with external 
angles $\theta_{i_1}$, $\ldots$, $\theta_{i_6}$ cyclically
is a Lorentz space with coordinates $(x, u, v, w)$ given by
\begin{equation*}
x = \sqrt{\Area T_0}, \quad
u = \sqrt{\Area T_1}, \quad  
v = \sqrt{\Area T_2}, \quad   
w = \sqrt{\Area T_3},
\end{equation*}
where $T_0$ is a triangle obtained by completing a convex hexagon 
by extending the edges $x_{i_2}$, $x_{i_4}$, $x_{i_6}$, 
and $T_1$, $T_2$, $T_3$ are yielding triangles adjacent to the edges 
$x_{i_1}$, $x_{i_3}$, $x_{i_5}$, respectively as in Fig.~\ref{Fig:10}.
\begin{figure}[ht]
\begin{center}
\setlength{\unitlength}{0.00083300in}%
\begin{picture}(3024,1824)(589,-2473)
\thicklines
\put(601,-2461){\line( 2, 3){1200}}
\put(1801,-661){\line( 1,-1){1800}}
\put(3601,-2461){\line(-1, 0){3000}}
\put(601,-2461){\line( 0, 1){  0}}
\put(1351,-1336){\line( 4,-1){1500}}
\put(1201,-1561){\line( 1,-3){300}}
\put(3001,-1861){\line(-2,-1){1200}}
\put(976,-2236){$T_1$}
\put(1201,-2386){$\theta_{i_2}$}
\put(2701,-2311){$T_2$}
\put(2926,-2086){$\theta_{i_4}$}
\put(2401,-1561){$\theta_{i_5}$}
\put(1801,-1936){$H$}
\put(1096,-1936){$\theta_{i_1}$}
\put(2273,-2393){$\theta_{i_3}$}
\put(1516,-1329){$\theta_{i_6}$}
\put(1779,-1186){$T_3$}
\end{picture}
\end{center}
\caption{}\label{Fig:10}
\end{figure}
Then $\Delta_{p,\theta}$ is a hyperbolic hexahedron 
in a hyperbola $\Area ^{-1}(1)$
in ${\cal E}_{p,\theta}$. By projecting it 
along the rays towards the origin into
the unit ball in $\{x=1\}$,  we get $\Delta_{ p,\theta}$ in 
the projective model
with $(u, v, w)$-coordinate.
Then $\Delta_{ p,\theta}$ lies in 
the first orthant  $\{u,v,w\, |\,u\geq 0, v\geq 0, w\geq
0\}$ with the vertex $(i_1i_2)(i_3i_4)(i_5i_6)$ at the origin and
the faces of $\Delta_{p,\theta}$ labeled
by $(i_1i_2)$, $(i_3i_4)$ and $(i_5i_6)$ lying 
in the planes $\{u=0\}$, $\{v=0\}$
and $\{w=0\}$, respectively.

The parameters $P$, $Q$, $R$ for the hexahedron $\Delta_{p,\theta} \in 
{\cal H}$ can be computed as follows.
Since $(P, 0,0)$ is the intersection of the plane $(i_6i_1)$ with
the $u$-axis, $(0,Q,0)$ the intersection of the plane $(i_2i_3)$
with $v$-axis, and  $(0,0,R)$ the intersection of the plane $(i_4i_5)$
with $w$-axis,
the points $(P,0,0)$, $(0,Q,0)$ and $(0,0,R)$ correspond 
to the degenerate hexagons with the sides
$x_{i_3}=x_{i_5}=x_{i_6}=0$,
$x_{i_1}=x_{i_2}=x_{i_5}=0$
and  $x_{i_1}=x_{i_3}=x_{i_4}=0$  respectively.

Let $a$, $b$, $c$ be the vertices of $T_0$ at which the
external angles are $\theta_{i_1}+\theta_{i_2}$, $\theta_{i_3}+\theta_{i_4}$,
$\theta_{i_5}+\theta_{i_6}$, respectively, and $c'$, $b'$, $a'$
the points on the (possibly extended) edges $ab$, $bc$, 
$ca$ such that the segments $cc'$, $aa'$, $bb'$
are parallel to the edges $x_{i_1}$, $x_{i_3}$, $x_{i_5}$ 
respectively (see Fig.\ref{Fig:11}).
\begin{figure}[ht]
\begin{center}
\setlength{\unitlength}{0.00083300in}%
\begin{picture}(3225,2271)(451,-2725)
\thicklines
\put(2601,-1461){\line(-2,-1){2000}}
\put(1801,-661){\line( 1,-3){600}}
\put(1037,-1820){\line( 4,-1){2564}}
\put(601,-2461){\line( 2, 3){1200}}
\put(1801,-661){\line( 1,-1){1800}}
\put(3601,-2461){\line(-1, 0){3000}}
\put(601,-2461){\line( 0, 1){  0}}
\put(1194,-1801){$\theta_{i_6}$}
\put(2513,-1719){$\theta_{i_4}$}
\put(2326,-2686){$c^\prime$}
\put(901,-1786){$b^\prime$}
\put(2626,-1411){$a^\prime$}
\put(1726,-586){$c$}
\put(451,-2611){$a$}
\put(3676,-2611){$b$}
\put(1711,-1043){$\theta_{i_1}$}
\put(3008,-2243){$\theta_{i_5}$}
\put(1133,-2363){$\theta_{i_3}$}
\put(2123,-2363){$\theta_{i_2}$}
\end{picture}
\end{center}
\caption{}\label{Fig:11}
\end{figure} 
Then by the similar procedure of presenting such degenerate hexagons in
the triangle $T_0$ as in the case $n=5$, we obtain

\begin{align*}
P^2 & = \frac{\Area \triangle ac'c}
             {\Area T_0} 
      = \frac{\text{length of } ac' }{\text{length of } ab}, \\
Q^2 & = \frac{\Area \triangle ba'a}
             {\Area T_0}
      = \frac{\text{length of } ba' }{\text{length of } bc }, \\
R^2 & = \frac{\Area \triangle cb'b}
             {\Area T_0}
      = \frac{\text{length of } cb' }{\text{length of } ca }.
\end{align*}

Thus we can compute $P$, $Q$, $R$ by just looking at the Fig.\ref{Fig:11}.
We note that the signs of $P-1$, $Q-1$, $R-1$ are the same as those of
$\theta_{i_5} +\theta_{i_6} + \theta_{i_1}-\pi$, 
$\theta_{i_1} +\theta_{i_2}+ \theta_{i_3}-\pi$, 
$\theta_{i_3} +\theta_{i_4}+ \theta_{i_5}-\pi$  respectively
(see Fig.\ref{Fig:12}).
\begin{figure}[ht]
\begin{center}
\setlength{\unitlength}{0.00083300in}%
\begingroup\makeatletter\ifx\SetFigFont\undefined%
\gdef\SetFigFont#1#2#3#4#5{%
  \reset@font\fontsize{#1}{#2pt}%
  \fontfamily{#3}\fontseries{#4}\fontshape{#5}%
  \selectfont}%
\fi\endgroup%
\begin{picture}(5487,2376)(526,-2125)
\thicklines
\put(601,-1561){\line( 1, 0){2100}}
\put(601,-1561){\line( 1, 2){900}}
\put(3301,-1561){\line( 1, 0){2700}}
\put(3301,-1561){\line( 1, 2){900}}
\put(1351,-61){\line( 3,-4){1125}}
\put(1351,-61){\line( 2,-5){596}}
\put(4051,-61){\line( 3,-4){1125}}
\put(4051,-61){\line( 6,-5){1788}}
%\multiput(3628,-921)(70.36364,-58.63636){12}{\makebox(6.6667,10.0000){\SetFigFont{10}{12}{\rmdefault}{\mddefault}{\updefault}.}}
%\multiput(1014,-756)(32.40000,-81.00000){11}{\makebox(6.6667,10.0000){\SetFigFont{10}{12}{\rmdefault}{\mddefault}{\updefault}.}}
\bezier{200}(1246,-264)(1321,-324)(1412,-287)
\bezier{200}(1443,127)(1550,-47)(1441,-186)
\bezier{200}(3956,-264)(4112,-321)(4224,-211)
\bezier{200}(4149,-189)(4249,-30)(4149,132)
\put(1201,-61){$c$}
\put(3901,-61){$c$}
\put(1876,-1786){$c^\prime$}
\put(2401,-1786){$b$}
\put(5101,-1786){$b$}
\put(5701,-1786){$c^\prime$}
\put(3226,-1786){$a$}
\put(526,-1786){$a$}
\put(1576,-61){$\theta_{i_5}+\theta_{i_6}$}
%%%% comment out the next line by inserting % to remove 
%\put(4351,-61){$\theta_{i_5}+\theta_{i_6}$}
%\put(1276,-1186){$x_{i_1}$}
%\put(4051,-1186){$x_{i_1}$}
\put(1276,-586){$\theta_{i_1}$}
\put(4051,-586){$\theta_{i_1}$}
\put(1201,-2086){(a) $P<1$}
\put(4276,-2086){(b) $P>1$}
\end{picture}
\end{center}
\caption{}\label{Fig:12}
\end{figure}

Now let $\Psi_p: \Theta_6 \to {\cal H}$ be the map sending $\theta \in \Theta_6$
to the hyperbolic hexahedron $\Delta_{p,\theta} \in {\cal H}$.
Then, as in the case $n=5$, the following result is 
proven in \cite{AharaYamada}, where we shall give 
a proof for reader's convenience.  

\begin{LemNum}\label{Lem:6fibration}
For a label $p=\langle i_1i_2\ldots i_6 \rangle$,
the inverse image 
of a point in ${\cal H}$ by $\Psi_p$  is in 
one-to-one correspondence with the upper half plane $\hyp^2$.
\end{LemNum}

\begin{pf}
Let $w$ be in $\hyp^2$.  Together with the point $0$ and $1$, it forms a
triangle which we denote by $T$.  
We add three edges in $T$.  
Let $X$ be the point  $P^2$  in the complex plane  
and add new edge $wX$.  
Likewise, let $Y$ be $1+(w-1)Q^2$ and
join $0$ and $Y$.  And last, let $Z$ be $w(1-R^2)$ and join $1$ and
$Z$.  By sliding the newly added edges like in the case $n=5$, we get
six angles $\theta_{i_1},\ldots,\theta_{i_6}$.  For the edge with
endpoints $X$, the direction to slide is toward $0$.  For $Y$ and $Z$,
it is $1$ and $w$ respectively (see Fig.~\ref{Fig:13}).  
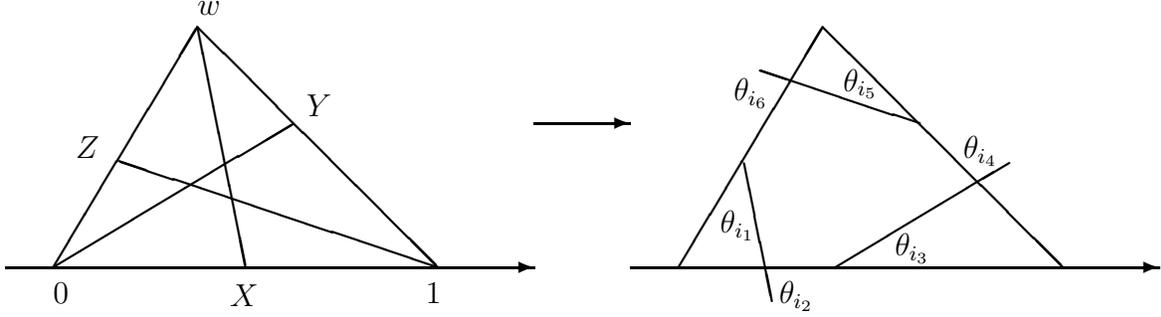
\begin{figure}[ht]
\begin{center}
\setlength{\unitlength}{0.00083300in}%
\begin{picture}(7224,1980)(589,-1834)
\put(5146,-518){$\theta_{i_6}$}
\thicklines
\put(4501,-1561){\vector( 1, 0){3300}}
\put(901,-1561){\line( 5, 3){1495}}
\put(3901,-661){\vector( 1, 0){600}}
\put(1801,-61){\line( 1,-5){300}}
\put(4801,-1561){\line( 3, 5){900}}
\put(5701,-61){\line( 1,-1){1500}}
\put(901,-1561){\line( 3, 5){900}}
\put(1801,-61){\line( 1,-1){1500}}
\put(3301,-1561){\line(-3, 1){1989}}
\put(5779,-1559){\line( 5, 3){1085}}
\put(6307,-663){\line(-3, 1){996}}
\put(5210,-908){\line( 1,-5){173}}
\put(1801, 14){$w$}
\put(901,-1786){$0$}
\put(6151,-1486){$\theta_{i_3}$}
\put(5431,-1786){$\theta_{i_2}$}
\put(3226,-1786){$1$}
\put(6578,-879){$\theta_{i_4}$}
\put(5829,-459){$\theta_{i_5}$}
\put(5063,-1344){$\theta_{i_1}$}
\put(2003,-1801){$X$}
\put(2483,-616){$Y$}
\put(1042,-879){$Z$}
\put(601,-1561){\vector( 1, 0){3300}}
\end{picture}
\caption{Map from $w$ to $\theta$}\label{Fig:13}
\end{center}
\end{figure}
Then by the argument above, we see that $\Psi_p((\theta_1, \ldots, \theta_6)) =
(P, Q, R)$. Since $w$ is arbitrary in $\hyp^2$, we get the conclusion.
\end{pf}

Consider now the map 
$\Psi_{\langle 123456 \rangle} \times \Psi_{\langle 214356\rangle}:\Theta_6 \to
{\cal H} \times  {\cal H}$ which sends $\theta$ to $(\Delta_{\langle 123456
\rangle,\theta}$, $\Delta_{\langle 214356 \rangle,\theta})$. 
Then these two hexahedra can be seen in Fig.~\ref{Fig:7}
stand on the first and third quadrants of the $uv$-plane. Notice that they have the
common edge $(12)(34)56$ along the $w$-axis. Then we claim the following.

\begin{LemNum}\label{Lem:6injective}
The map $\Psi_{\langle 123456 \rangle} \times \Psi_{\langle 214356 \rangle}$ is
injective.
\end{LemNum}

\begin{pf}
Let $((P_1,Q_1,R_1),(P_2,Q_2,R_2))$ be an element in the image of 
$\Psi_{\langle 123456 \rangle} \times 
\Psi_{\langle 214356 \rangle}$ and $\theta
= (\theta_1,\ldots,\theta_6)$ any element in the preimage
of $((P_1,Q_1,R_1),(P_2,Q_2,R_2))$. Set $w_1$ and
$w_2$ the points in $\xi_{\langle 123456 \rangle}^{-1}((P_1,Q_1,R_1))$ and
$\xi_{\langle 214356 \rangle}^{-1}((P_2,Q_2,R_2))$ corresponding to $\theta$
under the identification in Lemma~\ref{Lem:6fibration} 
respectively. 
Again, we see that the triangles $\triangle w_101$ and 
$\triangle w_201$ are
congruent so that $w_1$ and $w_2$ are identical
which we denote by $w$.

Let $X_i=P_i^2 $, $Y_i=1+ (w-1)Q_i^2 \in \hyp^2$ for $i=1,2$.
Then the triangles $\triangle 0X_1w$ and $\triangle 0wX_2$ are similar
since both have the external angles $\theta_{i_1}+\theta_{i_2}$,
$\pi-\theta_{i_2}$, $\pi-\theta_{i_1}$.
Similarly, the triangles $\triangle 10Y_1$ and $\triangle 1Y_20$ are similar
since both of which have the external angles $\theta_{i_3}+\theta_{i_4}$,
$\pi-\theta_{i_3}$, $\pi-\theta_{i_4}$.
(see Fig.~\ref{Fig:14}).
\begin{figure}[ht]
\begin{center}
\setlength{\unitlength}{0.00083300in}%
\begingroup\makeatletter\ifx\SetFigFont\undefined%
\gdef\SetFigFont#1#2#3#4#5{%
  \reset@font\fontsize{#1}{#2pt}%
  \fontfamily{#3}\fontseries{#4}\fontshape{#5}%
  \selectfont}%
\fi\endgroup%
\begin{picture}(4587,2898)(526,-3847)
\thicklines
\put(601,-3361){\line( 1, 0){1800}}
\put(601,-3361){\line( 1, 2){600}}
\put(3301,-3361){\line( 1, 0){1800}}
\put(3301,-3361){\line( 1, 2){600}}
\put(601,-3361){\line( 2, 1){1200}}
%\put(3256,-3103){\line( 3,-1){231}}
%\multiput(3533,-3344)(-2.10000,8.40000){21}{\makebox(6.6667,10.0000){\SetFigFont{7}{8.4}{\rmdefault}{\mddefault}{\updefault}.}}
\bezier{200}(3236,-3103)(3456,-3103)(3533,-3354)
\put(1201,-2161){\line( 1,-1){1200}}
\put(1201,-2161){\line( 1,-2){600}}
\put(5101,-3361){\line(-1, 1){2400}}
\put(3301,-3361){\line(-1, 4){600}}
\put(3901,-2161){\line( 2,-5){480}}
\put(3218,-3576){$0$}
\put(1426,-3811){(a)}
\put(4351,-3811){(b)}
\put(526,-3583){$0$}
\put(3601,-3231){$\theta_{i_4}$}
\put(3893,-2076){$w$}
\put(1193,-2076){$w$}
\put(1868,-2751){$Y_1$}
\put(1508,-3276){$\theta_{i_2}$}
\put(1718,-3576){$X_1$}
\put(2326,-3586){$1$}
\put(5026,-3586){$1$}
\put(4268,-3576){$X_2$}
\put(4118,-3276){$\theta_{i_1}$}
\put(3825,-2624){$\theta_{i_2}$}
\put(2476,-1111){$Y_2$}
\put(1141,-2542){$\theta_{i_1}$}
\put(1140,-3276){$\theta_{i_3}$}
\put(1748,-2999){$\theta_{i_4}$}
\put(2910,-1491){$\theta_{i_3}$}
\end{picture}
\end{center}
\caption{}\label{Fig:14}
\end{figure}
Thus we have
\[ P_1^2 : |w| = |w| : P_2^2,\]
and 
\[ 1 : |w-1|Q_1^2 = |w-1|Q_2^2 : 1. \]
Therefore $w$ is uniquely determined as the point on the upper half plane
where the two circles $|z|^2=(P_1P_2)^2$ and $|z-1|^2=1/(Q_1Q_2)^2$ meet.
And so we conclude that $\theta$ is unique.
\end{pf}

By the above lemma, we have proven the following theorem.

\begin{ThNum}
The map $\Phi_6$ from $\Theta_6$ to 
the space of Dehn fillings  ${\cal H}(\overline{X_6})$
defined by assigning to $\theta$ the hyperbolic cone manifold
$\overline{X_{6,\theta}}$ is injective.
\end{ThNum}

\end{document}